\documentclass[12pt]{article}%
\usepackage{amssymb}
\usepackage{amsfonts}
\usepackage{amsmath}
\usepackage{graphicx}%
\begin{document}

\begin{center}
\bigskip

\textbf{On distances and metrics in discrete ordered sets}
\end{center}

\bigskip

\begin{center}
Stephan Foldes

2017

\bigskip

\textbf{Abstract}
\end{center}

\textit{Discrete partially ordered sets can be turned into distance spaces in
several ways. The distance functions may or may not satisfy the triangle
inequality, and restrictions of the distance to finite chains may or may not
coincide with the natural, difference-of-height distance measured in a chain.
For semilattices, a semimodularity condition ensures the good behaviour of the
distances considered. }

\bigskip

Keywords: poset, lattice, semilattice, tree, semimodularity, chain condition,
height, distance, metric, triangle inequality

\bigskip

\bigskip

\textbf{1 \ Proximity in trees}

\bigskip

Degree of kinship between individuals has been considered relevant in ancient
and contemporary societies alike, in the normative context of laws of
inheritance, marriage prohibitions, rules against nepotism, and independence
of judges or jurors, to name a few examples. In Roman law, according to a
method now referred to as the "civil-law method", the degree of kinship
between two individuals, say Ego and Alter, is computed by determining their
nearest common ancestor X (which can be Ego or Alter if these two are in
direct line related), and then adding the number $h(E,X)$ of generations from
$X$ to Ego and the number $h(A,X)$ of generations from $X$ to Alter. According
to another ancient method, adopted in Europe in the Middle Ages and called the
"canon-law method" (see Bouchard for a historical account $\left[  B\right]  $
and Garner' legal dictionary $\left[  G\right]  $), the degree of kinship is
defined as the greater of the numbers $h(E,X)$ and $h(A,X)$. Our first
observation is that the triangle inequality is satisfied not only in the
computation of degree of kinship according to the civil-law method, but - less
obviously - in the computation according to the canon-law method as well.

\bigskip

\textbf{Proposition 1} \ \textit{Suppose that a set is partially ordered by a
discrete tree order} \textit{(i.e. all order-convex intervals }$[x,z]$%
\textit{\ are finite, and each pair }$x,y$\textit{\ of incomparable elements
has a least common upper bound }$x\vee y$ \textit{but has no common lower
bound). Then the distance function }$d(x,y)$ \textit{which assigns to elements
}$x,y$ \textit{the greater of }$(Card[x,x\vee y])-1$ \textit{and
}$(Card[y,x\vee y])-1$\textit{\ satisfies the triangle inquality }$d(x,y)\leq
d(x,v)+d(v,y).$ \ \ \ \ \ \ \ \ \ \ \ \ \ $\square$

\bigskip

In this note, abstracting from any possible applications or social context, we
shall formulate both the "civil-law" and "canon-law" methods of kinship degree
computation in the general abstract framework of partially ordered sets with a
connected Hasse diagram, examine the relationship between these and some other
distance functions, and address the question of validity of the triangle inequality.

\textbf{\bigskip}

\bigskip

\textbf{2 \ Distances in discrete partially ordered sets }

\bigskip

In the sequel a given partially ordered set, finite or infinite, shall be
called \textit{discrete}, if every maximal chain in every order-convex
interval $[x,y]$ is finite. This is a stronger condition than the requirement
that the order relation be generated as the transitive-reflexive closure of
its covering relation, which is a broader definition of discreteness adopted
for example in [FW]. However, discrete posets in the more restrictive sense
presently understood have the convenient property that the order induced on
any of their subsets is also discrete.

In a discrete poset, if two elements are comparable, say $x\leq y$, then by
the \textit{height} of $y$ above $x$, denoted indifferently by $h(x,y)$ or
$h(y,x)$ we mean the number that equals the least cardinality of a finite
maximal chain in $[x,y]$ minus $1.$

By a \textit{distance} function on a set $S$\ we mean a symmetric map $d$ from
$S^{2}$ to the non-negative reals for which $d(x,y)=0$ if and only if $x=y.$ A
distance function may or may not satisfy the \textit{triangle inequality}
\[
d(x,z)\leq d(x,y)+d(y,z)\text{
\ \ \ \ \ \ \ \ \ \ \ \ \ \ \ \ \ \ \ \ \ \ \ \ \ \ \ \ \ \ \ (1)}%
\]
while the term \textit{metric} is used for a distance function that
does.\textit{ }(Note that Deza and Laurent allow two distinct points to have
null distance, and the term distance is often used with the triangle
inequality being assumed to hold.) We shall also make use of a directed
(oriented) distance concept (see e.g. Chartrand, Johns, Tian and Winters
[CJTW], or Deza and Panteleeva [DP]), where the distance function is not
assumed to be symmetric, \textit{i.e.} where $d(x.y)$ may be different from
$d(y,x)$. Such a directed distance may also satisfy the triangle inequality.

The covering relation of any partial order defines a simple directed graph
with an arrow from element $x$ to element $y$ if and only if $x$ is covered by
$y$. Forgetting the orientation of the arrows, we obtain a simple undirected
graph called the poset's \textit{Hasse diagram}. If the Hasse diagram is
connected, then we call the poset \textit{connected}. Between any two elements
of a connected poset, we use the term \textit{zigzag distance} for their
graphic distance measured in the Hasse diagram of the poset. Zigzag distance
satisfies the triangle inequality $d(x,z)\leq d(x,y)+d(y,z)$ (this is so in
fact in non-discrete connected posets as well).

Recall that a poset has the \textit{upper} (respectively \textit{lower}%
)\textit{ filtering property }if any two elements have a common upper (lower)
bound. In a discrete poset with the upper (lower) filtering property, the
\textit{up--down }(respectively\textit{ down-up})\textit{ distance} of
elements $x$ and $y$ is defined as the smallest number of the form
$h(x,u)+h(y,u)$ (respectively of the form $h(u,x)+h(u,y)$ ), where $u$ is a
common upper (lower) bound of $x$ and $y.$ These notions are dual, trees and
other join semilattices have the upper filtering property, and lattices have
both filtering properties.

Obviously on any discrete chain the up-down, down-up and zigzag distance
functions coincide and yield what is conceivably the most natural notion of
distance on a chain. A distance function on a discrete poset is called
\textit{chain-compatible}, if its restriction to any maximal chain coincides
with this natural chain distance. This is a rather strong requirement, such
distance functions may not always exist:

\bigskip

\textbf{Proposition 2} \ \textit{For any discrete, connected partially ordered
set satisfying either one of the upper or lower filtering properties, the
following conditions are equivalent:}

\textit{(i) there is a chain-compatible distance function on the poset,}

\textit{(ii)\ the zigzag distance on the poset is chain-compatible,}

\textit{(iii) the poset satisfies the Jordan-Dedekind chain condition (in any
given interval} $\left[  x,y\right]  $ \textit{all maximal chains have the
same number of elements).}

\bigskip

\textbf{Proof} \ As each of the conditions \textit{(i) - (iii)} is self-dual,
we may suppose, without loss of generality, that the poset satisfies the upper
filtering condition.

Obviously condition \textit{(ii)} implies \textit{(i)}, and \textit{(i)}
implies \textit{(iii).} 

To show that \textit{(iii)} implies \textit{(ii)}, assume \textit{(iii)} and
suppose that there are elements $x<y$ for which the zigzag distance $d(x,y)$
is less than $h(x,y)$: this will lead to a contradiction. For each such pair
of elements $x<y$ there is a smallest positive integer $n=n(x,y)$, with the
property that there is a sequence of elements $x=x_{0},...,y=x_{n}$ , with
$x_{i}$ being comparable to $x_{i+1}$ for $0\leq i\leq n-1$, and such that
$h(x,y)>h(x_{0},x_{1})+...+h(x_{n-1},x_{n}).$ Choose $x<y$ so that $n=n(x,y)$
is minimal. Then $n\geq3,$ $x<x_{1}$ , $x_{1}>x_{2}$ and $x_{n-1}<y$. Let $u$
be a common upper bound of $x_{1}$ and $y.$ We must have, as $x_{n-1}<u$ and
$n$ is minimal,
\begin{align*}
h(x_{1},u) &  \leq h(x_{1},x_{2})+...+h(x_{n-1},u)=h(x_{1},x_{2}%
)+...+h(x_{n-1},y)+h(y,u)\\
h(x,y)+h(y,u) &  =h(x,u)=h(x,x_{1})+h(x_{1},u)\leq h(x,x_{1})+h(x_{1}%
,x_{2})+...+h(x_{n-1},y)+h(y,u)\\
h(x,y) &  \leq h(x,x_{1})+h(x_{1},x_{2})+...+h(x_{n-1},y)
\end{align*}

$\square$

\bigskip

A join semilattice is called \textit{semimodular} if for every pair of
distinct elements $x,y$, whenever there exists an element $z$ covered by both
$x$ and $y$, the join $x\vee y$ covers both $x$ and $y$. For lattices this
means just lattice semimodularity, but the extension obviously includes trees
as well. 

The following is easily seen to be true, by the usual argument for lattices:

\bigskip

\textbf{Proposition 3 \ }\textit{The Jordan-Dedekind chain condition is
satisfied in every discrete, semimodular join semilattice. \ }$\square$

\bigskip

Making use of this, again as in the case of lattices, we can see that for
discrete join semilattices, semimodularity is equivalent to the condition that
whenever elements $x,y$ have an element $z$\ as a common lower bound, we
should have $h(x,x\vee y)\leq h(z,y)$.\ 

\bigskip

Semimodularity is not necessary for the Jordan-Dedekind condition to hold.
Semimodularity can be characterized by the triangle inequality for the
\textit{non-symmetric semilattice distance} function $h(x,x\vee y)$:

\bigskip

\textbf{Proposition 4} \ \textit{A discrete join semilattice is semimodular if
and only if we have for all elements }$x,y,z\mathit{\ }$\textit{the
inequality}%
\[
h(x,x\vee y)+h(y,y\vee z)\geq h(x,x\vee z)\text{
\ \ \ \ \ \ \ \ \ \ \ \ \ \ \ \ \ \ \ \ \ \ \ \ \ \ \ \ \ \ (2)}%
\]

\bigskip

\textbf{Proof} \ Assume semimodularity. As $x$ is a common lowed bound of
$x\vee y$ and $x\vee z$ and $x\vee y\vee z$ \ is their join, we have%
\[
h(x,x\vee y)\geq h(x\vee z,x\vee y\vee z)
\]%
\[
h(y,y\vee z)\geq h(x\vee y,x\vee y\vee z)
\]
Therefore%
\[
h(x,x\vee y)+h(y,y\vee z)\geq\left[  h(x,h(x,x\vee y))-h(x\vee z,x\vee y\vee
z)\right]  +h(x\vee y,x\vee y\vee z)
\]
But the right hand side of this latter inequality equals $h(x,x\vee z).$

Conversely, if semimodularity fails, there are elements $x,y,z$ such that both
$x$ and $z$ cover $y$ but $x$ is not covered by $x\vee y$. Then
\[
h(x,x\vee y)=h(x,x)=0,\text{ \ \ }h(y,y\vee z)=h(y,z)=1,\text{ \ }h(x,x\vee
z)\geq2
\]
and (2) fails. \ \ \ \ \ \ $\square$

\bigskip

A further characterization of semimodularity can be given in terms of the
up-down distance. In any discrete poset, the up-down distance is always
greater than or equal to the zigzag distance, and if it satisfies the triangle
inequality, then it must be identical to the zigzag distance.

\bigskip

\textbf{Proposition 5} \ \textit{The following conditions are equivalent for
any discrete join semilattic}e $L$:

(i) $L$ \textit{is semimodular,}

(ii) \textit{the up-down distance on} $L$ \textit{satisfies the triangle
inequality,}

(iii)\textit{ the up-down distance on} $L$ \textit{coincides with the zigzag
distance.}

\bigskip

\textbf{Proof} \ First of all, if $L$ is not semimodular, then for some
elements $x,y,z$, the element $z$ is covered by both $x$ and $y$, but the join
$x\vee y$ does not cover $x$, i.e. $h(x,x\vee y)\geq2.$ Then $3\leq d(x,y)$
and\ the triangle inequality $d(x,y)\leq d(x,z)+d(z,y)=2$ fails for the
up-down distance.

Conversely, assume that $L$ is semimodular. If the triangle inequality failed
for the up-down distance, for some elements $x,y,z$ we would have%
\[
h(x,x\vee y)+h(y,x\vee y)+h(y,y\vee z)+h(z,y\vee z)<h(x,x\vee z)+h(z,x\vee z)
\]
But this is impossible, since by Proposition 4 we must have
\[
h(x,x\vee y)+h(y,y\vee z)\geq h(x,x\vee z)\text{
\ \ \ \ \ \ \ \ \ \ \ \ \ \ \ \ \ \ \ \ \ \ \ }%
\]

and%
\[
h(z,z\vee y)+h(y,y\vee x)\geq h(z,z\vee x)\text{
\ \ \ \ \ \ \ \ \ \ \ \ \ \ \ \ \ \ \ \ \ \ \ \ }\square
\]

\bigskip

On any discrete join semilattice, consider the "\textit{Chebyshev" distance}
function\textit{ }$d(x,y)=\max\left[  h(x,x\vee y),h(y,x\vee y)\right]  $
\ (The analogy - and overlap in the case of integer lattices - was pointed out
by Russ Woodroofe $\left[  \text{W}\right]  $.) Generally this distance need
not satisfy the triangle inequality. However, it does satisfy the triangle
inequality in a large class of semilattices, including trees (where it
corresponds to the canon-law method of determining degree of kinship). Note
that the Chebyshev distance, like the up-down distance, is always less than or
equal to the zigzag distance. 

\bigskip

\textbf{Proposition 6 \ }\textit{On any discrete, semimodular join
semilattice, the Chebyshev distance satisfies the triangle inequality.}

\bigskip

\textbf{Proof \ }Assume that the triangle inequality fails in some
semilattice, denote the Chebyshev distance by $d,$ and let $x,y,z$ be elements
such that $d(x,y)+d(y,z)<d(x,z).$ Let $a,b,c,d$ and $f,e$ denote the heights
$h(x,x\vee y),$ $h(y,x\vee y),$ $h(y,y\vee z),$ $h(z,y\vee z)$ and $h(x,x\vee
z),$ $h(z,x\vee z),$ respectively, in that order. Without loss of generality
$f\geq e$, and then $f$ must be (strictly) greater than each one of the
numbers $a+c,$ $a+d,$ $b+c,$ $b+d.$ Denote by $g,h,i$ the height of $x\vee
y\vee z$ above $x\vee y,$ $x\vee z,$ $y\vee z$, respectively. By the
Jordan-Dedekind condition, $f+h=a+g$. From this and from $f>a+c$ it follows
that
\[
a+c+h<a+g
\]
which implies $c<g.$ This contradicts semimodularity because $c=h(y,x\vee z)$
and $g=h\left[  x\vee y,(x\vee y)\vee(y\vee z)\right]  .$
\ \ \ \ \ \ \ \ \ \ \ \ \ \ $\square$

\bigskip

In contrast to to the equivalence of \textit{(i)} and \textit{(ii)} in
Proposition 5, semimodularity is only sufficient but not necessary for the
triangle inequality to hold for the Chebyshev distance in a discrete join
semilattice, as the example of the five element non-modular lattice shows. \ 

\bigskip

Finally, for any real $p\geq1$ consider the $l_{p}$ \textit{distance
function}, denoted $d_{p}$ on any discrete join semilattice, given by
\[
d_{p}(x,y)=\left[  h(x,x\vee y)^{p}+h(y,x\vee y)^{p}\right]  ^{1/p}\text{
\ \ \ \ \ \ \ \ \ \ \ \ \ \ \ \ (3)}%
\]
Obviously $d_{1}$ is the up-down distance, and - as expected - the Chebyshev
distance is the limit of the $l_{p}$\ distances as $p$ tends to infinity:%
\[
\underset{p\longrightarrow\infty}{\lim}d_{p}(x,y)=\max\left[  h(x,x\vee
y)+h(y,x\vee y)\right]
\]
In fact, again as expected, the $l_{p}$\ distance on any discrete semimodular
join semilattice \ (including all discrete semimodular lattices and trees)
satisfies the triangle inequality. In contrast to the Chebyshev distance,
semimodularity is characterized by the triangle inequality for any of the
$l_{p}$ distances on a discrete join semilattice, generalizing the equivalence
of \textit{(i)} and \textit{(ii)} in Proposition 5:

\bigskip

\textbf{Proposition 7} \ \ \textit{Let }$p\geq1$. \textit{A discrete join
semilattice }$L$\textit{\ is semimodular if and only if the }$l_{p}$\textit{
distance function }(3)\textit{ on }$L$\textit{\ satisfies the triangle
inequality.}

\bigskip

\textbf{Proof} \ Assume semimodularity. Proposition 4 allows to deduce the
triangle inequality from Minkowski's inequality (on which the triangle
inequality is based in classical $l_{p}$ spaces). In fact we only need the
following specialized two-dimensional case of Miskowski's inequality: if
$a_{1},a_{2},b_{1},b_{2}$ are non-negative real numbers and $1\leq p<\infty$,
then
\[
(a_{1}^{p}+a_{2}^{p})^{1/p}+(b_{1}^{p}+b_{2}^{p})^{1/p}\geq\left[
(a_{1}+b_{1})^{p}+(a_{2}+b_{2})^{p}\right]  ^{1/p}\text{
\ \ \ \ \ \ \ \ \ \ \ \ \ \ \ (4)}%
\]
To establish the triangle inequality for the $l_{p}$ distance $d_{p}$ in $L$
as defined by (3), we need to show that for all semilattice elements $x,y,z$%
\[
d_{p}(x,y)+d_{p}(y,z)\geq d_{p}(x,z)\text{
\ \ \ \ \ \ \ \ \ \ \ \ \ \ \ \ \ \ \ \ \ \ \ \ \ \ \ (5)}%
\]
Letting \ $a_{1}=h(x,x\vee y)$, \ $a_{2}=h(y,x\vee y)$, $\ b_{1}=h(y,y\vee
z)$, \ $b_{2}=h(z,y\vee z)$, the left had side of (5) is equal to the left
hand side of (4), while the right hand side of (4) is
\[
\left\{  \left[  h(x,x\vee y)+h(y,y\vee z)\right]  ^{p}+\left[  h(z,y\vee
z)+h(y,x\vee y)\right]  ^{p}\right\}  ^{1/p}\text{ \ \ \ \ \ (6)}%
\]
Now by Proposition 4
\[
\left[  h(x,x\vee y)+h(y,y\vee z)\right]  ^{p}\geq h(x,x\vee z)^{p}%
\]%
\[
\left[  h(z,y\vee z)+h(y,x\vee y)\right]  ^{p}\geq h(z,z\vee x)^{p}%
\]
and thus (6) is at least $d_{p}(x,y),$ completing the proof of (5).

Conversely, if semimodularity fails, then for some elements $x,y$ covering an
element $z$, the join $x\vee y$ does not cover $x$ and thus%
\begin{align*}
h(x,x\vee y)^{p}  & \geq2^{p}\\
d_{p}(x,y)^{p}  & >2^{p}\\
d_{p}(x,y)  & >2
\end{align*}
but $d_{p}(x,z)=d_{p}(y,z)=1$ and therefore (5) fails. $\ \square$

\bigskip

\bigskip

\textbf{Acknowledgements.}

This work, undertaken while the author was at the Tampere University of
Technology in Finland, has been co-funded by Marie Curie Actions (European
Union), and supported by the National Development Agency (NDA) of Hungary and
the Hungarian Scientific Research Fund (OTKA, contract number 84593), within a
project hosted by the University of Miskolc, Department of Analysis.

The author wishes to thank S\'{a}ndor Radeleczki and Russ Woodroofe for
valuable discussions.

\bigskip

\bigskip

\textbf{References}

\bigskip

[B] C.B. Bouchard, Consanguinity and noble marriages in the tenth and eleventh
centuries, \textit{Speculum} 56 (2) 268-287 (1981)

\bigskip

[CJTW] G. Chartrand, G.L. Johns, Songlin Tian, S.J. Winters, Directed distance
in digraphs: Centers and medians, \textit{J. of Graph Theory} 17 (4) 509--521 (1993)

\bigskip

[DL] M. Deza, M.\ Laurent, \textit{Geometry of Cuts and Metrics}, Springer 1997

\bigskip

[DP] M. Deza, E. Panteleeva, Quasi-semi-metrics, Oriented Multi-cuts and
Related Polyhedra, \textit{European J. of Combinatorics} 21 (6) 777-795 (2000)

\bigskip

[FW] \ S. Foldes, R. Woodroofe, \textit{Antichain cutsets of strongly
connected posets}, \textit{Order} 30 (2013), no. 2, 351-361

\bigskip

[G] \ B.A. Garner, \textit{A Dictionary of Modern Legal Usage}, Second
Edition, Oxford University Press (2001)

\bigskip

[W] R. Woodroofe, \textit{personal communication }(2013)

\end{document}